\documentstyle{article}

\def\C{{\bf C}}
\def\R{{\bf R}}
\def\F{{\bf F}}
\def\Z{{\bf Z}}

\def\Q{{\bf Q}}
\def\x{{\bf x}}
\def\y{{\bf y}}

\def\pf{\noindent{\bf Proof:\ }}
\def\qed{\hfill\framebox[2.5mm][t1]{\phantom{x}}}

\def\an{{\it Analogues:\/}\ }
\def\ans{{\it Analogues$^*$:\/}\ }

\def\halb{\hbox{$\frac{1}{2}$}}

\def\hexa{{\cal C}_6} 
\def\da{{\downarrow}}

\title{Self-dual Codes over the Kleinian Four Group}

\author{Gerald H\"ohn \\ \vspace{-1mm}
\small Mathematisches Institut der Universit\"at Freiburg \\ \vspace{-1mm}
\small Eckerstra{\ss}e 1, 79104 Freiburg, Germany \\ \vspace{-1mm}
\small E-mail: gerald@mathematik.uni-freiburg.de \\ \vspace{-1mm}
\small Homepage: http://baby.mathematik.uni-freiburg.de/  }

\date{26.~Mai 2000}

\begin{document}

\bibliographystyle{amsalpha} 

\newtheorem{theorem}{Theorem}
\newtheorem{lemma}{Lemma}

\renewcommand{\baselinestretch}{1.2}

\renewcommand{\thefootnote}{\fnsymbol{footnote}}

\maketitle

\section{Introduction}

In this work, we describe a new and natural fourth step in the series of analogies 
known to exist between binary codes, 
lattices and vertex operator algebras (see for example~\cite{CoSl,Ho-dr}).

Linear codes over the finite field ${\bf F}_4$ are studied in many papers
(cf.~\cite{MOSW-f4,CS-f3f4,Sloane-analogs,Sloane-complex,LP-no24,ChSl-aut18,Hu-ext1,Hu-ext2}), 
but a developed theory
for codes over the Kleinian four-group $K\cong \Z_2\times\Z_2$ is missing.
It turns out that there is a similar rich theory as one has for
binary linear codes. Parts of the results are known from some different 
viewpoints, but the use of Kleinian codes seems most natural.

We will prove all the results in terms of a theory for Kleinian codes, since this
leads to a theory of its own right, although one can deduce most theorems from 
the corresponding results for self-dual vertex operator algebras or lattices
or binary codes. To emphasize this relation,
we will give after every theorem a list of references of the
{\it analogue}\footnote{
An additional asterisk indicates that the theorem can be obtained from the analogues
theorems for binary codes, lattices or vertex operator algebras by the
relations described in the final section.}
theorems for binary codes (B), lattices (L) and vertex operator algebras (V).

The second section contains the main definitions and first results.
The next section describes the classification of odd and even self-dual
codes up to length $8$. In the fourth section, we study extremal codes.
This are codes with the largest possible minimal weight. The fifth section is about
designs for the space $K^n$.
Section six deals with lexicographic constructions. 

In the final section, we explain the relation and discuss some of 
the analogies with self-dual binary codes, lattices and vertex operator 
algebras in more detail. 
Self-dual Kleinian codes of length $n$ can be identified with self-dual vertex
operator superalgebras of rank $4n$ containing a vertex operator algebra of type 
$V_{D_4}^{\otimes n}$. From this viewpoint, Kleinian codes
are a special case of codes over a $3$-dimensional topological quantum 
field theory.

\medskip

Our motivation behind the introduction of Kleinian codes was to have an 
additional testbed besides binary codes and lattices for the understanding 
of vertex algebras. Kleinian codes have already found applications as
quantum codes and some of the results have been extended to and sharpened
for codes of larger length.

\section{Definitions and basic results}

Denote the elements of the Kleinian four group
$K\cong\Z_2\times\Z_2$ by $0$, $a$, $b$ and~$c$, where $0$ 
is the neutral element. The automorphism group of $K$ is $S_3$, the permutation
group of the three nonzero elements $a$, $b$ and~$c$. A {\it code\/} $C$ over $K$
of length~$n$ is a subset of the words of length~$n$ over the alphabet~$K$,
i.e.~consists of vectors $\x=(x_1,\ldots,x_n)$, $x_i\in K$, the {\it codewords\/}
of $C$. The {\it weight\/} ${\rm wt}(\x)$ of a codeword $\x$ is the number of 
nonzero $x_i$. The {\it minimal weight\/} of $C$ is defined by
$$ d=\min\{{\rm wt}(\x)\mid \x\in C,\,\ \x\not=0\}.$$
The code $C$ is called {\it linear\/} if $C$ is a subgroup of the abelian
group $K^n\cong\Z_2^{2n}$. A linear code has $4^k$ elements with $k\in\halb\Z$
and we denote $k$ the {\it dimension\/} of the code.
All codes in this article are assumed to be linear.
A code of length $n$, dimension $k$ and minimal weight $d$ is shortly denoted
as a $[n,k,d]$- or $[n,k]$-code. Let now $C$ be a $[n,k]$-code.

An important part of the structure which makes the theory of Kleinian codes 
interesting is the {\it scalar product\/} $(\,.\,,\,.\,)\,:\,K^n\times K^n
\longrightarrow \F_2$, $(\x,\y)=\sum_{i=1}^n x_i\cdot y_i$, where the symmetric
bilinear dot product $\cdot:K\times K\longrightarrow \F_2$ is defined by 
$a\cdot b=b \cdot a = 1$,
$a\cdot c=c \cdot a = 1$, $b\cdot c=c \cdot b = 1$ and zero otherwise.
The {\it dual code\/} $C^{\bot}$ is defined by
 $$C^{\bot}=\{\x\in K^n\mid (\x,\y)=0\ \hbox{for all\ }\y\in C\}$$ and has type
$[n,n-k]$.

We call $C$ {\it self-orthogonal\/} if $C\subset C^{\bot}$ and
{\it self-dual\/} if $C^{\bot}= C$.

The direct sum $C\oplus D$ of a $[n,k]$-code $C$ and a $[m,l]$-code $D$ is 
the direct product subgroup of $K^n\oplus K^m$ and has type $[n+m,k+l]$.
If $C$ can be written in a nontrivial way as a direct sum, $C$ is called 
{\it decomposable,\/} otherwise {\it indecomposable.} Obviously
$(C\oplus D)^{\bot}=C^{\bot}\oplus D^{\bot}$. Every code $C$ is after a 
renumbering of the positions a direct sum of indecomposable codes. 
The isomorphism classes of the components are uniquely determined up to 
permutation.

The {\it (Hamming) weight enumerator\/} of $C$ is the degree $n$ polynomial
$$W_C(u,v)=\sum_{i=0}^nA_i\,u^{n-i}v^i\qquad \hbox{with\ }
 A_i=\#\{\x\in C\mid {\rm wt}(\x)=i\}.$$
The {\it complete weight enumerator\/} is the polynomial
$${\rm cwe}_C(p,q,r,s)=\sum_{i,j,k,l}\, A_{i,j,k,l}\,p^iq^jr^ks^l,$$
where $A_{i,j,k,l}$ is the number of code words in $C$ containing at $i$,
$j$, $k$ resp.~$l$ of the $n$ positions the element $0$, $a$, $b$ resp.~$c$.
There is the obvious relation $W_{C}(u,v)={\rm cwe}_C(u,v,v,v)$.
Finally define for a natural number $g$ the {\it poly- or $g$-weight 
enumerator\/} $W^g_C$ as a polynomial in $2^g$ variables $t_{\nu}$ indexed
by~${\nu}\in \F_2^g$:
$$W^g_C=\sum_{\x^1,\,\ldots,\,\x^g\in C}\prod_{i=1}^n 
      t_{({\rm wt}(x^1_i),\ldots,{\rm wt}(x^g_i))} $$
and similar the {\it complete $g$-weight enumerator\/} ${\rm cwe}^g_C$ as a
polynomial in $4^g$ variables $s_{\nu}$ where~${\nu}\in K^g$.

The code $C$ is called {\it even\/} if the weights of all codewords are 
divisible by $2$. Note, that a code spanned by an orthogonal system of vectors
of even weight is itself even.

The automorphisms of the abelian group $K^n$ which are also 
isometries for the metric $d(\x,\y)={\rm wt}(\x-\y)$ on $K^n$
form the semidirect product $G=S_3^n{:}S_n$ consisting
of the permutation of the positions together with a permutation of the
symbols $a$, $b$ and~$c$ at each position. The {\it automorphism group\/} of
$C$ is the subgroup of $G$ sending $C$ to itself:
$$ {\rm Aut}(C)=\{g\in S_3^n{:}S_n \mid g\,C=C\}.$$
Two codes $C$ and $D$ are called to be equivalent if there is a $g\in G$ with
$g\, C=D$. The number of distinct codes equivalent to $C$ is
$$\frac{6^n\cdot n!}{|{\rm Aut}(C)|}. $$
Equivalent codes have the same (poly-) weight enumerator, but not necessarily the same
complete (poly-) weight enumerator. If $C$ is self-orthogonal, self-dual or even, so
it is the equivalent code.

Since $K$ is isomorphic to the additive group of the field $\F_4$, we can 
interpret every code over $K$ as a code over $\F_4$. Every code linear 
as $\F_4$-code is linear as Kleinian code, but not conversely. If $C$ is a self-dual 
Type IV $\F_4$-code for the hermitian scalar product of $\F_4^n$,
then it is also a even self-dual Kleinian code (cf.~\cite{MOSW-f4}). Perfect 
Kleinian codes are the same as perfect $\F_4$-codes, the only perfect 
Kleinian codes which exist
are $1$ error correcting codes~\cite{Ti-perfect}. 

\noindent{\bf Examples of Kleinian codes:}

- The $[1,\frac{1}{2},1]$-code $\gamma_1=\{(0),(a)\}$: $|{\rm Aut}(\gamma_1)|=2$,
$W_{\gamma_1}(u,v)=u+v$.

- The $[2,1,2]$-code $\epsilon_2=\{(00),(aa),(bb),(cc)\}$: 
$|{\rm Aut}(\epsilon_2)|=12$, $W_{\epsilon_2}(u,v)=u^2+3\,v^2$.

- The $[6,3,4]$-Hexacode ${\cal C}_6$ spanned by 
$$\{(a0a0bb),(a0bba0),(bba0a0),(00aaaa),(aa00aa),(b0b0ca)\}$$
as a Kleinian code. One has
$|{\rm Aut}({\cal C}_6)|=2^2\cdot 2\cdot 15\cdot 18=2160$,
$W_{{\cal C}_6}(u,v)=u^6+45\, u^2 v^4 + 18\, v^6$. 

- The Hamming code ${\cal H}_m$, $m\geq 2$ of type $[(4^m-1)/3,(4^m-1)/3-m,3]$
and the extended Hamming code $\overline{\cal H}_m$, $m\geq 2$ of type 
$[(4^m-1)/3+1,(4^m-1)/3-m,4]$. 

All examples are linear; the first three codes are self-dual;
$\epsilon_2$ and ${\cal C}_6\cong \overline{\cal H}_2$ and 
$\overline{\cal H}_m$ are even; besides $\gamma_1$, they are equivalent to codes
over ${\bf F}_4$; the code ${\cal H}_m$ is perfect.

\noindent{\bf Basic results:}

The Hamming weight enumerators of $C$ and its dual are related by the following
equation.
\begin{theorem}[generalized Mac-Williams identity (cf.~\cite{Del-four})]\label{mac-ham}
$$W_{C^{\bot}}(u,v)=\frac{1}{|C|}W_C(u+3v,u-v).$$
\end{theorem}
\ans B:~cf.~\cite{MacSl}; L:~cf.~\cite{Serre}, Ch.~VII, Prop.~16; 
V:~\cite{Zhu-dr,Ho-dr}.

\pf 
For a function $f$ on $K^n$ with values 
in a ring $R$ we define its transformation $g:K^n\longrightarrow R$ by
$g(\x)=\sum_{\y\in K^n}f(\y)\cdot (-1)^{(\y,\x)}$. One has the following
identity:
\begin{equation}\label{transform}
\frac{1}{|C|}\sum_{\x\in C}g(\x)=\sum_{\y\in C^{\bot}} f(\y).
\end{equation}\nopagebreak
Proof of~(\ref{transform}): $\sum_{\x\in C}g(\x)=\sum_{\y\in K^n}\sum_{\x\in C}
f(\y)(-1)^{(\x,\y)}=|C|\cdot \sum_{\y\in C^{\bot}} f(\y) + 
\sum_{\y\not\in C^{\bot}} f(\y)\cdot \sum_{\x\in C} (-1)^{(\x,\y)}$.
We have to show that the second sum vanishes. To this end, choose for
given $\y\in K^n\setminus C^{\bot}$ a $\x'\in C$ with $(\x',\y)\not=0$, 
i.e.,~$(-1)^{(\x',\y)}=-1$. We get 
$s=\sum_{\x\in C}(-1)^{(\x,\y)}=\sum_{\x\in C}(-1)^{(\x,\y)+(\x',\y)}=-s$,
which implies $s=0$ and proves~(\ref{transform}).

Now let $f(\y)=u^{n-{\rm wt}(\y)}v^{{\rm wt}(\y)}$. We obtain for its
transformation 
\begin{eqnarray*}
g(\x)& = &\sum_{\y\in K^n} f(\y)\cdot (-1)^{(\x,\y)} \\
     & = &\sum_{y_1,\ldots,y_n\in K} u^{n-{\rm wt}(y_1)-\cdots-{\rm wt}(y_n)}
v^{{\rm wt}(y_1)+\cdots+{\rm wt}(y_n)}(-1)^{x_1\cdot y_1+\cdots + x_n\cdot y_n}\\
 & =& \prod_{i=1}^n\left(\sum_{z\in K} u^{1-{\rm wt}(z)}v^{{\rm wt}(z)}(-1)^{x_i\cdot z}\right)\\
 & = & (u+3v)^{n-{\rm wt}(\x)}(u-v)^{{\rm wt}(\x)}.
\end{eqnarray*}
Applying (\ref{transform}) we get for the weight enumerator of $C^{\bot}$:
$$
W_{C^{\bot}}(u,v) = \sum_{\y\in C^{\bot}}f(\y) 
 =  \frac{1}{|C|}\sum_{\x\in C} (u+3v)^{n-{\rm wt}(\x)}(u-v)^{{\rm wt}(\x)}
 =  \frac{1}{|C|} W_C(u+3v,u-v). $$
\qed

For the other types of weight enumerators we stay only the results, 
the proofs are similar.
\begin{theorem}[Mac-Williams identity for complete weight enumerators] 
$${\rm cwe}_{C^{\bot}}(p,q,r,s)=
\frac{1}{|C|}{\rm cwe}_C(p+q+r+s,\,p+q-r-s,\,p-q+r-s,\,p-q-r+s).$$
\end{theorem}

From Theorem~\ref{mac-ham}, we get the following descriptions of the weight 
enumerators of self-dual codes:
\begin{theorem}\label{gleason-odd}
Let $C$ be a self-dual $[n,n/2]$-code. Then, the weight enumerator
$W_C(u,v)$ is a weighted homogeneous polynomial of weight $n$ in $u+v$ and
$v(u-v)$, or equivalently in the weight enumerators 
of $\gamma_1$ and $\epsilon_2$.
\end{theorem} 
\ans B:~\cite{Gleason}; L:~cf.~\cite{CoSl}; V:~\cite{Ho-dr}, Ch.~2.

\pf From Theorem~\ref{mac-ham}, we see that $W_C$ is invariant under
the group $H\cong\Z_2$ generated by the substitution 
$\frac{1}{2}{{1\ \phantom{-}3 \choose  1 \ -1}}$. The ring of 
invariants has Molien series $1/\big((1-\lambda)(1-\lambda^2)\big)$.
(This is the generating function for the multiplicities of the trivial
$H$-representation in the symmetric powers of the defining two
dimensional representation of $H$.) The polynomials $u+v$ and
$v(u-v)$ or equivalently $W_{\gamma_1}$ and $W_{\epsilon_2}$ 
are algebraically independent and generate freely the ring of all invariants.
\qed 

\begin{theorem}\label{gleason-even}
Let $C$ be an even self-dual $[n,n/2]$-code. Then, the weight enumerator
$W_C(u,v)$ is a weighted homogeneous polynomial of weight $n$ in 
$u^2+3v^2$ and $v^2(u^2-v^2)^2$, or equivalently in
the weight enumerators of $\epsilon_2$ and ${\cal C}_6$.
\end{theorem} 
\ans B:~\cite{Gleason}; L:~cf.~\cite{CoSl};
V: see~\cite{Go-mero} and~\cite{Ho-dr}, Ch.~2.

\nopagebreak
\pf This follows from the corresponding result for even self-dual
codes over $\F_4$ as proven for example in~\cite{MOSW-f4}, Th.~13: 
The group generated by
$S=\frac{1}{2}{{1\ \phantom{-}3 \choose  1 \ -1}}$ and
$T={{1\ \phantom{-}0 \choose  0 \ -1}}$ has order $12$ and
the Molien series of the corresponding ring of invariants is 
$1/\big((1-\lambda^2)(1-\lambda^6)\big)$. \phantom{neuezeile} \qed

\section{Classification of self-dual codes}

Let $\delta_n$ be the code consisting of all codewords containing only $0$'s and an 
even number of $a$'s. This is the even subcode of $\gamma_1^n$.
One has ${\rm dim}(\delta_n)=(n-1)/2$, coset representatives
of $\delta_n^{\bot}/\delta_n$ are given by $(0^n)$, $(a,0^{n-1})$, $(b^n)$ and 
$(c,b^{n-1})$ and its automorphism group consists for $n\geq 2$ of the permutation 
of the positions together with possible interchanging $b$ and $c$ at every 
position, i.e.,~${\rm Aut}(\delta_n)=S_2^{n}{:}S_n$.

The next theorem describes self-orthogonal codes spanned by vectors of small weight.
\begin{theorem}\label{smallgen}
Minimal weight $1$ subcodes of a self-orthogonal code $C$ can be split off: 
$C\cong D\oplus\gamma_1^l$, with minimal weight of $D$ larger then $1$.
Self-orthogonal codes generated by weight-$2$-vectors are equivalent to direct 
sums of $\delta_l$, $l\geq 2$, and $\epsilon_2$.
\end{theorem}
\ans First part: B, L: easy to see; 
V: cf.~\cite{Go-mero}, \cite{Ho-dr}, Th.~2.2.8.
Second Part: B: \cite{PlSl}, Th.~6.5; L: cf.~\cite{CoSl}, Ch.~4; V: Cartan, Killing,
\cite{infinite}, \cite{FreZhu}.
 
\pf For the first statement, note that a weight-$1$-codeword 
is equivalent to $(0,\ldots,0,a)$. Then $C$=$C'\oplus \gamma_1$, where $C'$ is
the orthogonal complement in $C$ of the $\gamma_1$ spanned by $(0,\ldots,0,a)$.

For the proof of the second statement, decompose first the code
generated by the weight-$2$-codewords into the direct sum of its
indecomposable even components and fix one of them. We have two
possibilities:

Case i) There are two weight-$2$-codewords containing different nonzero entries  
at the same position. \newline
In this case the component is equivalent to a code containing the two
codewords $(aa0\ldots0)$ and $(bb0\ldots0)$. They generate a $\epsilon_2$
subcode and, since $\epsilon_2$ is self-dual, this is the whole component.
(The other possible pairs of weight-$2$-codewords are not
orthogonal.)

Case ii) The component is equivalent to a code whose weight-$2$-codewords have at all 
positions the value $0$ or $a$. \newline
Inductively, one sees that the component is equivalent to a $\delta_l$, 
$l\geq 2$. A possible set of generators is given by $(aa0\ldots0)$,
$(0aa0\ldots)$, $\ldots$, $(0\ldots0aa)$. \qed

Let $\bar C$ the subcode of $C$ generated by the weight $1$ and $2$ codewords. We 
can describe $C$ by its {\it gluecode\/} $\Lambda\subset \bar C^{\bot}/\bar C$.
The automorphism group of $C$ is given by ${\rm Aut}(C)=G_0. G_1.G_2$,
where $G_0$ are the ``inner automorphisms'' of $\bar C$, i.e.~those which
are fixing the components of $\bar C$ and the cosets $\Lambda/\bar C$, $G_1$ are
the automorphisms of $C$ fixing the components of $\bar C$ modulo $G_0$ and 
$G_2$ is the induced permutation group on the components of $\bar C$.

Denote by $M(n)$ resp.~$M_e(n)$ the number of distinct (but maybe equivalent) 
self-dual resp.~even self-dual Kleinian codes of length $n$.
\begin{theorem}[Massformula]\label{mas-tot}
The mass constants are given by
$$M(n)=\prod_{i=1}^{n}(2^i+1)=\sum_{C}\frac{6^n\cdot n!}{|{\rm Aut}(C)|} $$
where the sum is over equivalence classes of self-dual codes and
$$M_e(n)=\prod_{i=0}^{n-1}(2^i+1)=\sum_{C}\frac{6^n\cdot n!}{|{\rm Aut}(C)|} $$
where the sum is over equivalence classes of even self-dual codes and $n$ is 
even.
\end{theorem}
\an B: cf.~\cite{PlSl}; L:~\cite{min}; V: unknown. 

\pf First, we prove the formula for $M(n)$. Let $M(n,k)$ be the 
number of self-orthogonal codes of dimension $k$ and length $n$. There are
$|(C^{\bot}\setminus C)/C|=4^{n-2k}-1$ different extensions of a self-orthogonal 
$[n,k]$-code $C$ to a self-orthogonal $[n,k+\frac{1}{2}]$-code $D\supset C$
by choosing one extra vector $x\in C^{\bot}$. Every  self-orthogonal
$[n,k+\frac{1}{2}]$-code $D$ arises from $|D\setminus\{0\}|=4^{k+1/2}-1$ 
different codes $C$. So we get the recursion
$$M(n,k+\frac{1}{2})=M(n,k)\cdot\frac{4^{n-2k}-1}{4^{k+1/2}-1}.$$
Together with $M(n,0)=1$ we obtain
$$ M(n)=M(n,n/2)=\prod_{i=0}^{n-1}\frac{4^{n-i}-1}{2^{i+1}-1}=
\prod_{i=1}^{n}(2^i+1).$$ 
The second expression for $M(n)$ describes the decomposition of all self-dual 
codes into orbits under the action of $S_3^n{:}S_n$.

To get the mass formula for $M_e(n)$, define in a similar way as before
$M_e(n,k)$ as the number of even self-orthogonal codes of dimension $k$ and 
length $n$. The dual code $C^{\bot}$ of a even self-orthogonal 
$[n,k]$-code $C$ contains $\frac{1}{2}(4^{n-k}+(-2)^n)$ vectors of even weight
as one can see from Theorem~\ref{mac-ham}. All vectors in a coset 
$C^{\bot}/C$ have the same weight modulo $2$. So we get in a similar way as above the 
recursion
$$M_e(n,k+\frac{1}{2})=M_e(n,k)\cdot\frac{\frac{1}{2}(4^{n-2k}+2^{n-2k})-1}
{4^{k+1/2}-1}.$$
Starting from $M_e(n,0)=1$ we obtain
$$ M_e(n)=M_e(n,n/2)=\prod_{i=0}^{n-1}
\frac{2^{2n-2i-1}+2^{n-i-1}-1}{2^{i+1}-1}=\prod_{i=0}^{n-1}(2^i+1)$$ 
and again we can express the total number as a sum over the different 
equivalence classes of codes.
\qed

For the weighted sum of the Hamming weight enumerators one has
\begin{theorem}[Massformula for Hamming weight enumerators]\label{mas-ham}
$$ \sum_{C}\frac{6^n\cdot n!}{|{\rm Aut}(C)|}W_C(u,v)=
M(n)\cdot(1+2^{n})^{-1}\cdot\left[2^n u^n+(u+3v)^n\right] $$
where the sum is over equivalence classes of self-dual codes.
$$\sum_{C}\frac{6^n\cdot n!}{|{\rm Aut}(C)|}W_C(u,v)= 
M_e(n)\cdot(1+2^{n-1})^{-1}\cdot\left[2^{n-1} u^n+\frac{1}{2}\left\{(u+3v)^n
+(u-3v)^n\right\}\right] $$
where the sum is over equivalence classes of even self-dual codes.
\end{theorem}
\an B:~\cite{PlSl}; L:~\cite{sie}; V: unknown.

\pf Let $\x$ be a nonzero vector (of even weight) of length $n$. Similar
as in the proof of Theorem~\ref{mas-tot} one gets for the number of
(even) self-dual codes containing $\x$ the expression
$$\prod_{i=1}^{n-1}(2^i+1)\quad 
\hbox{or\ \ }\prod_{i=0}^{n-2}(2^i+1)\hbox{\ for even codes.}$$
From this and Theorem~\ref{mas-tot} one obtains the result by summing
$$ u^{n-{\rm wt}(\x)} v^{{\rm wt}(\x)} $$ over all pairs $(\x,C)$, where $C$
is a (even) self-dual code with $\x\in C$, and expanding the resulting
sum in two different ways. \qed

We remark, that the average Hamming weight enumerator for even self-dual
Kleinian codes is the same as for even formal self-dual $\F_4$-codes 
(\cite{MOSW-f4}, Th.~24) although the mass constants are different.

We call a self-dual code {\it primitive\/}, if no $\gamma_1$ subcode can be split off.
A primitive code $C$ is the first one in the chain $C$, $C\oplus\gamma_1$, 
$\ldots$
\begin{theorem}[Relation between even and odd self-dual codes]\label{relation}
There is a $1:1$-correspondence between isomorphism classes of
pairs $(C,\delta_k)$, where $C$ is an
even self-dual code of even length $n$ and $\delta_k$  a subcode (a defined above) 
inside $C$
(together with the choice of a class $[x]$ in $\delta_k^{\bot}/\delta_k$
of minimal weight $1$, i.e., if $k=1$ we must choose
$x\in K\setminus\{0\}$) and 
isomorphism classes of self-dual codes $D$ of length $n-k$. 

The code $D$ is primitive if and only if the subcode $\delta_k$ is maximal,
i.e.~not contained in a $\delta_{k+1}$ subcode (with corresponding gluevectors
$[x]$).
\end{theorem}
\ans B:~\cite{CoPless}; L:~\cite{CoSl-ul23}; 
V:~\cite{Ho-dr}, Ch.~3, and~\cite{Ho-svoa}.

\pf We describe the map from self-dual codes $D$ of length $n-k$ to even self-dual
codes of even length $n$. Denote by $\delta_k^0=\delta_k$, $\delta_k^1$,
$\delta_k^2$ and $\delta_k^3$ the four cosets of $\delta_k$ inside 
$\delta_k^{\bot}$ such that $(a0^{k-1})\in \delta_k^1$.

If $D$ is even, let $C=D\oplus(\delta_k^0\cup\delta_k^2)$. 
Otherwise we have the decomposition 
$D_0^{\bot}=D_0\cup D_1\cup D_2\cup D_3$
of the orthogonal complement of the even subcode $D_0$ of $D=D_0\cup D_1$ 
into four $D_0$ cosets. Define
$$C=D_0\oplus \delta_k^0\cup D_1\oplus \delta_k^1\cup
    D_2\oplus \delta_k^2\cup D_3\oplus \delta_k^3. $$
Note that for $k=1$ the three cosets $\delta_k^1$, $\delta_k^2$ and $\delta_k^3$
are all equivalent under ${\rm Aut}(\delta_1)=S_3$.
It is then easy to check that this map describes the claimed 
$1:1$-correspondence.
\qed 

We call $D$ a {\it child\/} of the parent code $C$. From Theorem~\ref{relation}, we get 
the following description of the primitive children of an even 
self-dual code $C$ of length $n$: Take a position and choose $x\in\{a,b,c\}$ 
(up to the action of ${\rm Aut}(C)$), this gives a self-dual code $D$ of length $n-1$.

- If the position is not in the support of the subcode $\bar C$ generated by the 
weight-$2$-codewords, the code $D$ is primitive.

- If the position is in an $\delta_l$, $l\geq 2$, component of $\bar C$ we have two 
cases: If $x\not=a$ then $D$ is again maximal, if $x=a$ the primitive child is
obtained by deleting the remaining $l-1$ positions of $\delta_l$ from $D$.

- If the position is in a $\epsilon_2$ component, the primitive child is
obtained by deleting the second position of $\epsilon_2$ from $D$.

\smallskip

\begin{figure}\caption{The neighbourhood graph 
for $n=2$, $4$ and $6$}\label{neighbgraph}
\vspace{2mm}

\setlength{\unitlength}{1mm}
\begin{picture}(140,70)
\put(0,67){$n=2$:}

\put(9,59){$\epsilon_2$}\put(10,60){\circle{8}}

\thicklines
\put(19,60){\line(-1,0){5} } \put(20,59){$\gamma_1^2$}

\thinlines

\put(0,50){$n=4$:}

\put(8,39){$\delta_4^+$}\put(10,40){\circle{8}}
\put(9,19){$\epsilon_2^2$}\put(10,20){\circle{8}}

\put(10,36){\line(0,-1){12} } \put(11,28){$\delta_2^2$}

\thicklines
\put(14,40){\line(1,0){5} }  \put(17.5,47){$\gamma_1^4$}
\put(13,43){\line(1,1){4} }  \put(20,39){$\delta_3\gamma_1$}
\put(14,20){\line(1,0){5} }  \put(20,19){$\epsilon_2\gamma_1^2$}
\thinlines

\put(40,0)
{\begin{picture}(100,70)
\put(0,67){$n=6$:}

\put(48,59){$\delta_6^+$}\put(50,60){\circle{8}}

\put(16.5,39){$\delta_4^+\epsilon_2$}\put(20,40){\circle{8}}

\put(78,39){${\delta_3^2}^+$}\put(80,40){\circle{8}}

\put(19,19){$\epsilon_2^3$}\put(20,20){\circle{8}}

\put(78,19){${\delta_2^3}^+$}\put(80,20){\circle{8}}

\put(49,-1){${\cal C}_6$}\put(50,0){\circle{8}}

\put(53,57){\line(3,-2){23} } \put(26.5,50){$\delta_4\delta_2$}
\put(23,43){\line(3,2){23} } \put(65,50){$\delta_3^2$}
\put(24,40){\line(1,0){52} } \put(21,28){$\epsilon_2\delta_2^2$}
\put(24,39){\line(3,-1){52} }  \put(50.5,31){$\delta_2^2\delta_2$}
\put(20,36){\line(0,-1){12} } \put(76,28){$\delta_2^2$}
\put(80,36){\line(0,-1){12} } \put(48,41.5){$\delta_3\delta_2$}
\put(24,20){\line(1,0){52} } \put(48,21.5){$\delta_2^3$}
\put(53,3){\line(3,2){23} } \put(58,10){${\cal O}_6$}

\thicklines

\put(54,60){\line(1,0){5} }  \put(37,59){$\gamma_1^6$}
\put(46,60){\line(-1,0){5} } \put(60,59){$\delta_5 \gamma_1$}

\put(17,43){\line(-1,1){4} } \put(8,48.5){$\epsilon_2 \gamma_1^4$}
\put(16,40){\line(-1,0){5} } \put(4,39){$\delta_4 \gamma_1^2$}
\put(17,37){\line(-1,-1){4} } \put(7,30.5){$\delta_3\epsilon_2\gamma_1$}

\put(84,40){\line(1,0){5} }  \put(87.5,47){$\delta_3\gamma_1^3$}
\put(83,43){\line(1,1){4} }  \put(90,39){$\delta_3\delta_2\gamma_1$}

\put(17,17){\line(-1,-1){4} } \put(8,9){$\epsilon_2^2\gamma_1^2$}

\put(84,20){\line(1,0){5} }  \put(89.5,19){$\delta_2^2\gamma_1^2$}
\put(83,17){\line(1,-1){4} } \put(86,9.5){$\delta_2^2\gamma_1$}
\put(80,16){\line(0,-1){5} } \put(79,8){$\delta_2$}

\put(46,0){\line(-1,0){5} } \put(34.5,-1){${\cal C}_5\gamma_1$}

\end{picture}}
\end{picture}

\end{figure}
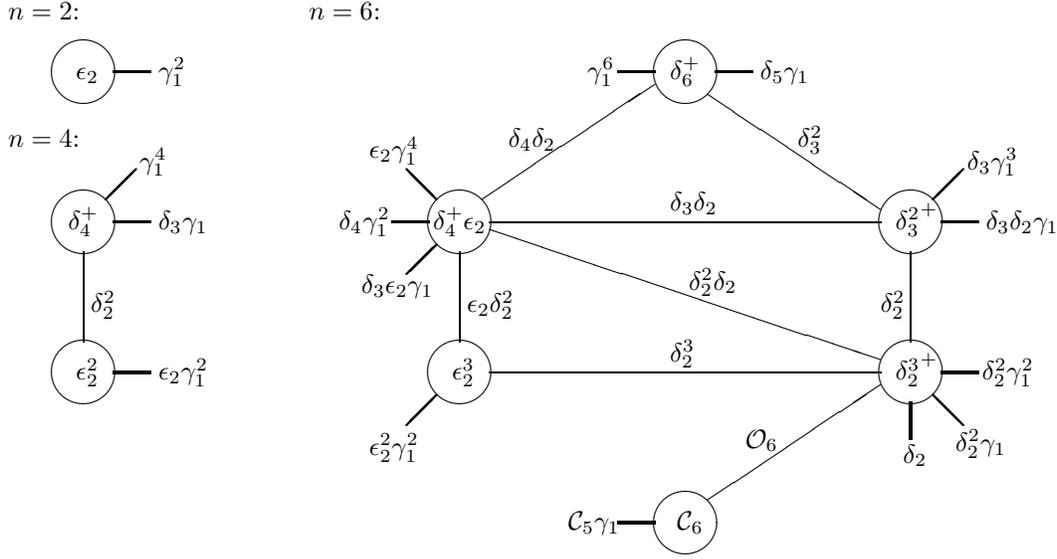

Every non even self-dual code $D=D_0\cup D_1$ of even length $n$ determines
the two {\it even\/} self-dual ``neighbours'' $D_0\cup D_2$ and 
$D_0\cup D_3$, where $D_0$, $D_1$, $D_2$ and $D_3$ are the four cosets
of $D_0$ in $D_0^{\perp}$ as above. We define for every even $n$ a
``neighbourhood graph'' by using the isomorphism classes of even self-dual
codes as vertices, the isomorphism classes of non even self-dual
codes as edges and ``neighbourhood'' as incidence relation. An edge corresponding
to a non primitive code $D=D'\oplus \gamma_1^l$, $l\geq 1$, is a loop for
the vertex corresponding to the even code determined from $D'$ through
Theorem~\ref{relation}. The edges starting on a vertex $C$ correspond to
the orbits of ${\rm Aut}(C)$ on the nonzero elements of $K^n/C$.
It is easy to see that the neighbourhood graph is connected for all $n$. 
For $n=2$, $4$ and $6$ the graph is shown in Figure~\ref{neighbgraph}.

\noindent\ans L:~\cite{Bo-thesis}; V: not determined.

\medskip

\begin{theorem}\label{class-even}
The even self-dual codes up to length $8$ (together with the subcode $\bar C$,
order of $G_1.G_2$, weight enumerator and number of children)
are given in Table~\ref{even-codes}.
\begin{table}\caption{Even self-dual codes up to length $8$\label{even-codes}}
\vspace{5mm}
$$\begin{array}{|rr|lr|rrrrr|rr|}
n & \hbox{No.} & \bar{C} & |G_1||G_2| & A_0 & A_2 & A_4 & A_6 & A_8 & n_1 & n_2 \\ \hline\hline
2 & 1 & \epsilon_2    & 1 & 1 & 3 & & & & 1 & 1\\ \hline
4 & 1 & \delta_4      & 1 & 1 & 6 & 9 & & & 2 & 1 \\
  & 2 & \epsilon_2^2  & 2 & 1 & 6 & 9 & & & 1 & 1 \\ \hline
6 & 1 & \delta_6           & 1  & 1 & 15& 15 & 33 & & 2 & 1 \\
  & 2 & \delta_4\epsilon_2 & 1  & 1 & 9 & 27 & 27 & & 3 & 2 \\
  & 3 & \delta_3^2         & 2  & 1 & 6 & 33 & 24 & & 2 & 1 \\
  & 4 & \epsilon_2^3       & 3! & 1 & 9 & 27 & 27 & & 1 & 1 \\
  & 5 & \delta_2^3         & 3! & 1 & 3 & 39 & 21 & & 2 & 1 \\
  & 6 & {\cal C}_6         &2160& 1 & 0 & 45 & 18 & & 1 & 0 \\ \hline
8 & 1 & \delta_8&1          & 1& 28& 70& 28&129& 2& 1 \\ 
& 2& \delta_6\epsilon_2&1   & 1& 18& 60& 78& 99& 3& 2 \\
& 3& \delta_5\delta_3&1     & 1& 13& 55&103& 84& 4& 2 \\
& 4& \delta_4^2&2           & 1& 12& 54&108& 81& 2& 1 \\
& 5& \delta_4\epsilon_2^2&2 & 1& 12& 54&108& 81& 3& 2 \\
& 6& \delta_4\delta_2^2&2   & 1&  8& 50&128& 69& 4& 2 \\
& 7& \delta_3^2\epsilon_2&2 & 1&  9& 51&123& 72& 3& 2 \\
& 8& \delta_3^2\delta_2&2   & 1&  7& 49&133& 66& 4& 2 \\
& 9& \delta_3\delta_2^2&2   & 1&  5& 47&143& 60& 7& 2 \\
& 10& \delta_3&120          & 1&  3& 45&153& 54& 3& 1 \\
& 11& \epsilon_2^4&24       & 1& 12& 54&108& 81& 1& 1 \\
& 12& \epsilon_2\delta_2^3&6& 1&  6& 48&138& 63& 3& 2 \\
& 13& \epsilon_2 & 2160     & 1&  3& 45&153& 54& 2& 1 \\
& 14& \delta_2^4&24         & 1&  4& 46&148& 57& 2& 1 \\
& 15& \delta_2^4&8          & 1&  4& 46&148& 57& 2& 1 \\
& 16& \delta_2^3&6          & 1&  3& 45&153& 54& 3& 1 \\
& 17& \delta_2^2&16         & 1&  2& 44&158& 51& 4& 1 \\
& 18& \delta_2&48           & 1&  1& 43&163& 48& 4& 1 \\
& 19& -&6\cdot1344          & 1&  0& 42&168& 45& 1& 0 \\
& 20& -&1152                & 1&  0& 42&168& 45& 1& 0 \\
& 21& -&336                 & 1&  0& 42&168& 45& 1& 0 \\ \hline
\end{array}$$
\end{table}
\end{theorem}
\ans B:~\cite{CoPless,CoPlSl-bin32rev}; L:~\cite{kneser,niemeier}; 
V: cf.~\cite{Go-mero}, for $c=24$ there is a conjectured list in
\cite{Schellekens1}.

\pf Use the list of doubly even self-dual binary codes of length 
$4n$~\cite{CoPless,CoPlSl-bin32rev} and the construction~A described in 
Section~\ref{KBLV} or use Theorem~\ref{smallgen} and 
classify the possibilities for $\bar{C}$ and the gluecodes $\Lambda\subset
\bar{C}^{\bot}/\bar{C}$ directly.
\qed

We checked the result additionally with the mass formula for the 
Hamming weight enumerator. 

\begin{theorem}\label{class-odd}
The non even self-dual codes up to length $6$ (together with the parent No.,
the subcode $\bar C$, order of $G_1.G_2$ and the weight enumerator) 
are given in Table~\ref{odd-codes}.
\begin{table}\caption{The non even self-dual codes up to length $6$\label{odd-codes}}
\vspace{5mm}
$$\begin{array}{|rr|r|lr|rrrrrrr|}
n & \hbox{No.} & \hbox{par.\ No.} &\bar{C} & |G_1||G_2| & A_0 & A_1 & A_2 & A_3 & A_4 & A_5 &  A_6  \\ \hline\hline
1 & 1 & 1 & \gamma_1            & 1 & 1 & 1 &&&&& \\ \hline
2 & 1 & 1 & \gamma_1^2     &     2! & 1 & 2 & 1 &&&&\\ \hline
3 & 1 & 1  & \gamma_1^3       &  3! & 1 & 3 & 3 & 1 & & &\\ 
  & 2 & 2 & \epsilon_2\gamma_1 &  1 & 1 & 1 & 3 & 3&&&\\ 
  & 3 & 1  & \delta_3          &  1 & 1 & 0 & 3 & 4&&&\\ \hline
4 & 1 & 1  & \gamma_1^4          & 4! & 1 & 4 & 6 & 4 & 1&& \\
  & 2 & 2  & \epsilon_2\gamma_1^2 &2! & 1 & 2 & 4 & 6 & 3&& \\ 
  & 3 & 3  & \delta_3\gamma_1     & 1 & 1 & 1 & 3 & 7 & 4&& \\ 
  & 4 & 5  & \delta_2^2          & 2! & 1 & 0 & 2 & 8 & 5&& \\  \hline
5 & 1& 1& \gamma_1^5&5!&            1 & 5 & 10&  10& 5& 1 &\\  
  & 2& 2& \epsilon_2\gamma_1^3 &3!& 1 & 3 &6 &10& 9& 3&\\
  & 3& 3& \delta_3\gamma_1^2&2!&    1 & 2 &4& 10& 11& 4&\\
  & 4& 5& \delta_2^2\gamma_1&2!&    1 & 1 &2& 10& 13& 5&\\
  & 5& 2& \delta_4\gamma_1&1&       1 & 1 &6& 6& 9& 9&\\ 
  & 6& 4& \epsilon_2^2\gamma_1&2!&  1 & 1 &6& 6& 9& 9&\\
  & 7& 1&  \delta_5&1&              1 & 0 &10&0&5&16&\\
  & 8& 2& \delta_3\epsilon_2&1&     1 & 0 &6&4&9&12&\\
  & 9& 3& \delta_3\delta_2&   1&    1 & 0 &4&6&11&10&\\
  &10& 5& \delta_2^2&2!&            1 & 0 &2&8&13&8&\\
  &11& 6& {\cal C}_5 & 120     &    1 & 0 &0&10&15&6&\\ \hline
6 & 1& 1&\gamma_1^6             & 6! &1& 6& 15& 20& 15& 6& 1\\
  & 2& 2&\epsilon_2\gamma_1^4 & 4!   &1& 4& 9& 16& 19& 12& 3\\
  & 3& 3&\delta_3\gamma_1^3      & 3!&1& 3& 6& 14& 21& 15& 4\\
  & 4& 6&\delta_2^2\gamma_1^2 & 2!^2 &1& 2& 3& 12& 23& 18& 5\\
  & 5& 4&\delta_4\gamma_1^2      &2! &1& 2& 7& 12& 15& 18& 9\\
  & 6& 5&\epsilon_2^2\gamma_1^2& 2!^2&1& 2& 7& 12& 15& 18& 9\\
  & 7& 3&\delta_5\gamma_1&          1&1& 1& 10& 10& 5& 21& 16\\
  & 8& 7&\delta_3\epsilon_2\gamma_1&1&1& 1& 6& 10& 13& 21& 12\\
  & 9& 8&\delta_3\delta_2\gamma_1&  1&1& 1& 4& 10& 17& 21& 10\\
  &10& 9&\delta_2^2\gamma_1&      2! &1& 1& 2& 10& 21& 21& 8\\
  &11&10&{\cal C}_5\gamma_1&       120&1& 1& 0& 10& 25& 21& 6\\
  &12& 6&\delta_4\delta_2& 1         &1&0&7&8&7&24&17\\
  &13& 8&\delta_3^2       &2!        &1&0&6&8&9&24&16\\
  &14& 9&\delta_3\delta_2& 1         &1&0&4&8&13&24&14\\
  &15&12&\epsilon_2\delta_2^2&2!     &1&0&5&8&11&24&15\\
  &16&14&\delta_2^3       &3!        &1&0&3&8&15&24&13\\
  &17&15&\delta_2^2\delta_2&2        &1&0&3&8&15&24&13\\
  &18&16&\delta_2^2       &2         &1&0&2&8&17&24&12\\
  &19&17&\delta_2         &8         &1&0&1&8&19&24&11\\
  &20&18&{\cal O}_6       &6\cdot 8  &1&0&0&8&21&24&10\\ \hline
\end{array}$$
\end{table}

\end{theorem}
\ans B:~\cite{Pless-small,PlSl}; L:~\cite{CoSl-ul23,Bo-l24}; 
V:~\cite{Ho-dr}, Ch.~3, and~\cite{Ho-svoa}. 

\pf Look at the list of even self-dual binary codes of length~$4n$~\cite{Pless-small,PlSl} 
or apply Theorem~\ref{relation} to Theorem~\ref{class-even}. \qed

Again we checked the result by the mass formula for the Hamming weight 
enumerator. 

{\it Remark:\/} There is one self-dual code of length $5$ without codewords
of weight $2$: The {\it shorter Hexacode\/} ${\cal C}_5$. There are two self-dual
codes of length $6$ without codewords of weight~$2$: The Hexacode 
${\cal C}_6$ (even) and the {\it odd Hexacode\/} $O_6$ (non even). 

The number of inequivalent (even) self-dual codes of small length~$n$ can
be read off from Table~\ref{anzahlcodes}.
\begin{table}\caption{Number of inequivalent (even) self-dual codes}\label{anzahlcodes} 
\vspace{5mm}
\begin{tabular}{l|cccccccccccccc}
Type & $1$ & $2$ & $3$ & $4$ & $5$ & $6$ & $7$ & $8$ & $9$ & $10$ & $11$ & $12$\\ \hline
even & - & $1$ & - & $2$ & - & $6$ & - & $21$ & - &  & - & $\geq 338$ \\ \hline
odd  & $1$  & $2$ & $3$ & $6$ & $11$ &  $26$ & $59$ & & & $\geq 392$ & $\geq 12143$ \\ \hline
\end{tabular}
\end{table}
The number of even codes up to length~$8$ are obtained from 
Theorem~\ref{class-even}, the number of odd codes up to length~$6$
from Theorem~\ref{class-odd} and~\ref{class-even} and for $n=7$ it follows
from the number of length~$7$ children of the even length~$8$ codes.
The lower estimates for larger~$n$ one obtains from the mass 
formula.

A complete classification up to $n=10$ seems possible, but no interesting
new structure is expected. 

\smallskip

All the self-dual Kleinian codes classified in this section
have a nontrivial automorphism group. In analogy to~\cite{OrPh,Bannai-rigid},
we expect that this holds only for small length $n$ and that rather
almost all self-dual and even self-dual codes have trivial automorphism group.
What are the smallest (even) self-dual codes with trivial automorphism groups
(cf.~\cite{Ba-trivial} for lattices)?

\section{Extremal codes}

In this section, we study self-dual Kleinian codes of type $[n,n/2,d]$ where
$d$ is as large as possible.
\smallskip
Let $m=[n/2]$. By Theorem~\ref{gleason-odd} the weight enumerator of a 
code $C$ can be written as
\begin{equation}\label{def-ai}
W_C(u,v)=\sum_{i=0}^{m}a_i\,(u+v)^{n-2i}(v(u-v))^i
\end{equation}
with unique integral numbers $a_i$.
There is a unique choice of the numbers $a_0$, $\ldots$, $a_m$ such that
the right hand side of~(\ref{def-ai}) equals
\begin{equation}\label{def-ext}
u^n+0\cdot u^{n-1}v+\cdots + 0\cdot u^{n-m}v^m + A_{m+1}\,u^{n-m-1}v^{m+1}
+ \cdots + A_n\, v^n.
\end{equation}
We call~(\ref{def-ext}) the {\it extremal weight enumerator\/} and a code
with this weight enumerator {\it extremal}.
So an extremal code has minimal weight $d\geq[n/2]+1$.

\begin{theorem}\label{bound-odd}
The minimal distance $d$ of a self-dual code $C$ of length $n$ satisfies
$$d\leq\left[\frac{n}{2}\right]+1.$$
\end{theorem}
\an B:~\cite{MaSl}; L:~\cite{Siegel-zeta}; V:~\cite{Ho-dr}, Cor.~5.3.3.

\pf The proof is parallel to~\cite{MaSl}, Cor.~3. In fact it can be 
considered as ``case~5''\footnote{``Case 4'' was defined in~\cite{MOSW-f4}.}
of that paper for the parameters $w=1$, $R=2$, $S=1$ and $\alpha=1$. It follows
also from the next theorem. \qed

Let $C_0$ be the even subcode of $C$ as in the proof of Theorem~\ref{relation}.
To study extremal codes in more detail, we need the definition of the
{\it shadow\/} $C'$ of $C$: 
We set $C'=C_0^{\bot}\setminus C$ if $C$ is not even and $C'=C$ otherwise. 

\begin{lemma}\label{shadow-lemma}
If the weight enumerator of $C$ is written as
$$ W_C(u,v)=P_C(W_{\gamma_1},W_{\epsilon_2})=
Q_C(W_{\gamma_1},W_{\epsilon_2}-W_{\gamma_1^2}) $$
with weighted homogeneous polynomials $P_C(x,y)$ and $Q_C(x,y)$, then for 
the shadow one has
$$ W_{C'}(u,v)=P_C(W_{\gamma'_1},W_{\epsilon_2})=
Q_C(W_{\gamma'_1},W_{\epsilon_2}-W_{{\gamma'_1}^2}).$$
\end{lemma}
\pf We show $W_{C'}(u,v)=\frac{1}{|C|} W_C(u+3v,(-1)(u-v))$ from which the
lemma follows.

If $C=C'$ is even this is Theorem~\ref{mac-ham}. Otherwise, we get from there
\begin{eqnarray*}
W_{C'}(u,v)& = & W_{C_0^{\bot}}(u,v)-W_C(u,v)\ =\ 
\frac{2}{|C|}W_{C_0}(u+3v,u-v)-W_C(u,v)\\
& = & \frac{1}{|C|}\left[W_{C}(u+3v,u-v)+W_{C}(u+3v,(-1)(u-v))\right]-W_C(u,v)\\
&=&\frac{1}{|C|}W_{C}(u+3v,(-1)(u-v)).\qquad \qed
\end{eqnarray*}  

\begin{theorem}\label{class-ext}
There are exactly five extremal codes: $\gamma_1$, $\epsilon_2$, $\delta_3^+$,
the shorter Hexacode ${\cal C}_5$ and the Hexacode ${\cal C}_6$.
\end{theorem}
\an B:~\cite{MaSl,Wa}; L:~\cite{CoOdSl}; V:~\cite{Ho-dr}, Th.~5.3.2.

For the corresponding extremal weight enumerators see Table~\ref{even-codes} 
and~\ref{odd-codes}.

\pf The existence and uniqueness of an extremal code for $n=1$, $2$, $3$,
$5$ and $6$ can directly be read off from Table~\ref{even-codes} 
and~\ref{odd-codes}.

The nonexistence for $n=4$ follows also from this tables, so we must
prove the nonexistence for $n>6$. We can assume $C$ is non even since for
an even code we will show (Theorem~\ref{bound-ext-even}) that for the minimal 
weight $d$ one has $d\leq 2[n/6]+2$. But from $d\geq [n/2]+1$, we get $n=2$
or $6$. Now we are using the shadow $C'$ of $C$.
From Lemma~\ref{shadow-lemma}, we get for its weight enumerator
for $n=7$, $8$, $\dots$, $11$:
$$
\begin{array}{l|ccccc}
n           & 7 & 8 & 9 & 10 & 11  \\ \hline
W_{C'}(u,v) & \phantom{\frac{|}{|}}\frac{7}{4}u^6v+\cdots & -\frac{13}{8}u^8+\cdots & 
-\frac{9}{4}u^8v+\cdots &\frac{23}{8}u^{10} +\cdots &\frac{33}{8}u^{10}v +\cdots  
\end{array}.
$$
Since $W_{C'}$ must have non negative integral coefficients, there exists
no extremal codes for $7\leq n\leq 11$. For $n\geq 12$, the
coefficient $A_{m+2}$ of $W_C(u,v)$ is always negative. We will sketch the
proof:

Let $m=[n/2]$ and replace $u$ by $1$. Expanding $(1+v)^{-n}$ in
powers of $\phi=\frac{v(1-v)}{(1+v)^2}$ one gets by
the B\"urmann Lagrange Theorem
\begin{equation}\label{expan}
(1+v)^{-n}=\sum_{k=0}^{m}b_k\,\phi^k+\sum_{k=m+1}^{\infty} b_k\,\phi^k
\end{equation}
with 
$$b_k=\left.\frac{1}{k!}\frac{d^{k-1}}{dv^k}\left[\frac{d(1+v)^{-n}}{dv}
\left(\frac{v}{\phi}\right)^k\right]\right|_{v=0}.$$
Comparing expansion~(\ref{expan}) with~(\ref{def-ai}) and (\ref{def-ext})
yields $b_k=a_k$ for $k=0$, $\ldots$, $m$. Furthermore, $A_{m+1}=-b_{m+1}$,
$A_{m+2}=-b_{m+2}+3(m+1)b_{m+1}-n$. Now one estimates with the saddle-point 
method $b_{m+1}$ and $b_{m+2}$ and shows that $A_{m+2}<0$ for $m$ large enough. 
The smaller $n$ are checked by a direct computation.
\qed

{\it Remarks:\/} Similar as in~\cite{CoSl-f1,CoSl-JN,CoSl-JN-korrekt} one can 
refine the bound of Theorem~\ref{bound-odd} to obtain 
\hbox{$d\leq 2[n/5]+{\rm O}(1)$} by using the shadow code.

For the difference $D_C(u,v)=W_{C_2}(u,v)-W_{C_3}(u,v)$ one has the result
$$
D_C(u,v)\in 
\cases{\Q[W_{\epsilon},W_{{\cal C}_6}], & if $n$ is even, \cr
    v(u^2-v^2)\Q[W_{\epsilon},W_{{\cal C}_6}], & if $n$ is odd.}  
$$
This result can be used as in~\cite{CoSl-f1,CoSl-JN,CoSl-JN-korrekt} 
to discuss for small $n$ the ``weakly'' extremal codes meeting the 
stronger bound for $d$. As an example, for $n=5$ we obtain 
$D_C(u,v)=c\cdot v(u^2-v^2)(u^2+3 v^2)$.  

\smallskip

Instead of looking for codes with large minimal weight, one can ask the
same question for the shadow itself.
For self-dual codes with shadows of large minimal weight one gets similar 
results as recently described by N.~Elkies and the author:
\begin{theorem}\label{shadow-strict} 
The minimal weight $h$ of the shadow $C'$ of a self-dual code $C$ of length $n$
satisfies $h\leq n$, with equality if and only if $C\cong \gamma_1^n$.
\end{theorem}
\ans B:~\cite{El-shadow}; L:~\cite{El-Z}; V:~\cite{Ho-shadow}, Th.~1.

\pf Clearly $h\leq n$.
By Lemma~\ref{shadow-lemma}, the weight enumerator of $C'$ is a polynomial
$P_C(W_{\gamma'_1},W_{\epsilon_2})$ in the weight enumerators of $\gamma_1'$
and $\epsilon_2$, i.e.~$W_{C'}(u,v)$ is a homogeneous polynomial of weight
$n$ in $2v$ and $u^2+3v^2$. So $h=n$ implies $W_{C'}(u,v)=(2v)^n$; but then
$W_C(u,v)=(u+v)^n$ and $C\cong\gamma_1^n$. \qed

\begin{theorem}\label{shadow-lax}
Let $C$ be a self-dual code of length $n$ without words of weight $1$.
Then one has
\begin{description}
\item [i)] $C$ hat at least $(n/2)(5-n)$ codewords of weight $2$.
\item [ii)] The equality holds if and only if $h(C')=n-2$.
\item [iii)] In this case the number of codewords of weight $n-2$ in the shadow
is $2^{n-3}\cdot n$.
\end{description}
\end{theorem}
\ans B, L:~\cite{El-shadow}; V:~\cite{Ho-shadow}, Th.~2.

\pf Assume first $h(C')\geq n-2$. In the same way as in the proof of 
Theorem~\ref{shadow-strict} we see that $P_C(x,y)$  is a linear combination
of $x^n$ and $x^{n-2}y$ and we obtain
\begin{eqnarray}\label{weight-lax}
W_C(u,v) & = & (u+v)^n - \frac{n}{2}(u+v)^{n-2}\left( (u+v)^2-(u^2+3v^2)\right) \\
 & = & u^n+ 0\cdot u^{n-1}v + \frac{n}{2}(5-n)u^{n-2}v^2+\cdots .
\end{eqnarray} 
This proves one direction of ii). 

Conversely, we can assume $n<6$, so the
weight enumerator of $C$ can be written as
\begin{eqnarray*} 
W_C(u,v) &\! =\! & (u+v)^n - \frac{n}{2}(u+v)^{n\!-\!2}\left(2uv-2v^2\right)
 + \frac{A_2-(n/2)(5\!-\!n)}{4} (u+v)^{n\!-\!4}\left(2uv-2v^2\right)^2.
\end{eqnarray*}
From Lemma~\ref{shadow-lemma}, we get $A_2-(n/2)(5-n)\geq 0$ since 
$W_{C'}(u,v)$ has nonnegative coefficients, and we have i) and the
converse of of ii).

Finally, Part iii) follows also from (\ref{weight-lax}) and
Lemma~\ref{shadow-lemma}:
\begin{eqnarray*}
W_{C'}(u,v) & = & (2v)^n -\frac{n}{2}(2v)^{n-2}\left((2v)^2-(u^2+3v^2)\right) \\
 & =& 2^{n-3} n\cdot u^2v^{n-2} + \left(2^n-n\,2^{n-3}\right)v^n. \qquad \qed
\end{eqnarray*} 

There are exactly four such codes meeting the bound $h(C')=n-2$, namely 
$\epsilon_2$, $\delta_3^+$, $(\delta_2^2)^+$ and ${\cal C}_5$.

\medskip

For even codes there are similar definitions and results. 
The following result was proven for ${\bf F}_4$-codes, but since its proof uses only
Theorem~\ref{gleason-even} it is also true for Kleinian codes.
\begin{theorem}[see~\cite{MOSW-f4}]\label{bound-ext-even}
The minimal distance $d$ of an even self-dual code $C$ of length $n$ satisfies
$$d\leq2\left[\frac{n}{6}\right]+2.$$
\end{theorem}
\an B:~\cite{MaOdSl}; L:~\cite{MaOdSl}; V:~\cite{Ho-dr}, Section~5.2.

{\it Remark:\/} The analogous bound for doubly-even binary codes 
has recently been improved in~\cite{KrLi,KrLi2} for large lengths.

\smallskip

An even self-dual code matching this bound is called {\it extremal\/}.
The corresponding weight enumerator is called the {\it extremal weight 
enumerator\/} of length $n$. A table of extremal weight enumerators 
was given in~\cite{MOSW-f4}, Table~1.

Again from the ${\bf F}_4$ case, the next result follows.
\begin{theorem}[see~\cite{MOSW-f4}]\label{bound-even}
There are no extremal even codes of length $n\geq 136$.
\end{theorem}
\an B:~\cite{MaOdSl}; L:~\cite{MaOdSl}; 
V:~no known bound, cf.~\cite{Ho-dr}, Section~5.2.

Examples of extremal ${\bf F}_4$-codes are known for
$n=2$ ($\epsilon_2$), $4$ ($\epsilon_2^2$), $6$ ($\hexa$), $8$ ($3$ codes), $10$,
$14$, $\ldots$, $22$, $28$ and $30$ (see~\cite{CS-f3f4}). 
They are also examples of extremal even Kleinian codes. 

There is no extremal ${\bf F}_4$-code of length~12.
But there is an extremal even Kleinian code of this length with
generator matrix
$$\left(
\begin{array}{ll}
\tt aaaaaa & \tt 000000 \\
\tt bbbbbb & \tt 000000 \\
\tt 000000 & \tt aaaaaa \\
\tt 000000 & \tt bbbbbb \\
\tt a0bab0 & \tt aaaa00 \\
\tt abccba & \tt bbbb00 \\
\tt caca00 & \tt a0aaa0 \\
\tt cca0a0 & \tt b0bbb0 \\
\tt ccbaab & \tt a00aaa \\
\tt bccbaa & \tt b00bbb \\
\tt caabcb & \tt aa00aa \\
\tt b0baa0 & \tt bb00bb \\
\end{array}
\right)$$
and weight enumerator
$W_C(u,v)=u^{12}+396\,u^6v^6+1485\,u^4v^8+1980\,u^2v^{10}+234\,v^{12}$.

Besides the question of the existence of a projective plane of order ten
and of a doubly even code of type $[72,36,16]$, a $[24,12,10]$ self-dual
${\bf F}_4$-code was most wanted. 
After the first question, also the third question has 
found a negative answer~\cite{LP-no24}. Since Kleinian codes are combinatorial
more natural than ${\bf F}_4$-codes, we ask if there is 
an even self-dual Kleinian code of type $[24,12,10]$.
This is the smallest open case for extremal even Kleinian codes.

\bigskip

Good even and doubly even self-dual binary codes meeting the Gilbert-Varshamov 
bound exist, as was shown by using the mass formula for the Hamming weight 
enumerator~\cite{wo}. A similar result holds for lattices (see~\cite{MiHu}, Ch.~II). 
We expect the same for self-dual and even self-dual Kleinian codes.

\section{Constant weight codes and generalized $t$-designs}

Let $X_k$ be the fiber over $k$ of the weight map ${\rm wt}:K^n\longrightarrow
\{0,1,\ldots,n\}$. We can write it as the (not two point) homogenous space
$X_k=G/H=S_3^n{:}S_n/(S_2^k{:}S_k\times S_3^{n-k}{:}S_{n-k})$. 
The $H$-module structure of the function space $L_2(X_k)$ for general 
alphabets instead of $K$ has been studied in~\cite{dunkl}.
The space $X_k$ carries the structure of a symmetric association scheme, called
the nonbinary Johnson scheme (cf.~\cite{TaAaGo}) as follows: A pair 
$(\x,\y)\in X_k\times X_k$ belongs to the relation $R_{r,s}$, with $r$, $s\in 
\{0,1,\ldots, k\}$, $r\leq s$, if $r=\#\{i\mid x_i=y_i\not=0\}$ and
$s=\#\{i\mid x_i\not=0,\ y_i\not=0\}$. This structures allow one to use the 
usual association scheme methods to study subsets $Y\subset X_k$ (cf.~\cite{DeLe}).
\footnote{I like to thank C.~Bachoc for mentioning the
references~\cite{dunkl,TaAaGo,DeLe,Del-four,Ba-harm} to me.}

\medskip

Here, we use the definition of a generalized $t$-designs as 
in~\cite{Del-four}: An element $\x\in K^n$ is said to be covered
by an element $\y\in K^n$ if each nonzero component $x_i$ of $\x$ is
equal to the corresponding component $y_i$ of $\y$. 
A {\it generalized $t-(n,k,\mu)$ design
(of type $3$)\/} is a nonempty subset $Y\subset X_k$ such that any element of $X_t$
is covered by exactly $\mu$ elements from $Y$. 
For $t=2$, this definition is identical with the notion of a 
{\it group divisible incomplete block design\/} with $n$ groups of $3$ elements,
blocksize $k$ and $\lambda_1=0$, $\lambda_2=\mu$ introduced in~\cite{BoNa}. 

As an example, the three codewords of weight $2$ in $\epsilon_2$ form
a generalized $1$-$(2,2,1)$ design. The next result describes a method to
obtain generalized $2$-designs.

\begin{theorem}\label{design}
Let $C$ be an extremal even code of length $n=6k$. Then, 
the codewords of $C$ of fixed non-zero weight form a generalized $2$-design.
\end{theorem}
\an B:~\cite{AsMa}; L:~\cite{Ve-design}; V: unknown.

\pf This follows from Th.~5.3.\ in~\cite{Del-four}, a generalization
of the Assmus and Mattson theorem: By Theorem~\ref{bound-ext-even},
there are at most $\frac{1}{2}\big(n-(2(n/6)+2)\big)+1=2(n/6)$ nonzero 
weights in such a code. Note, that our scalar product on $K$
defines a required identification map $\chi_{(.)}:K\longrightarrow {\rm Hom}(K,{\bf C}^*)$.
\qed

The result applies in particular to the unique extremal even code of 
length~$6$, the Hexacode ${\cal C}_6$ and the extremal even code of 
length~$12$ given in the last section.
The generalized $2$-$(6,4,2)$ and  $2$-$(6,6,2)$ designs formed by the vectors 
of the Hexacode of weight $4$ and $6$ are unique. 

\medskip

In this case, the design property can also be obtained 
from the following result about ${\rm Aut}({\cal C}_6)$:

\begin{theorem}\label{authexa-wt2}
The automorphism group of the Hexacode acts transitively
on the weight~$2$ vectors in $K^6$.
\end{theorem}
\an B:~\cite{Mathieu1861,Carmwitt}; L:~\cite{GoSe,HuSl}; V: unknown.

\pf By computing the double cosets 
$(S_2^k{:}S_k\times S_3^{6-k}{:}S_{6-k})\setminus S_3^6/{\rm Aut}({\cal C}_6)$ for 
\hbox{$k=0$, $1$, $\ldots$, $6$}, we get the orbit decomposition of
$K^6$ under ${\rm Aut}({\cal C}_6)$ as shown in Table~\ref{hexaorbits}.
\begin{table}\caption{Orbits of ${\rm Aut}({\cal C}_6)$ in $K^6$}\label{hexaorbits}
$$\begin{array}{|c|c|r|r|r|}
\hbox{weight $k$} & \hbox{name}  & \hbox{Size} & \hbox{Distance to ${\cal C}_6$} &
\hbox{nearest codeword(s)}\\ \hline \hline
0 & A_0 &   1 & 0 & A_0 \\ \hline
1 & A_1 &  18 & 1 & A_0\\ \hline
2 & A_2 & 135 & 2 & A_0, 2\times A_4\\ \hline
3 & A_3 & 180 & 1 & A_4\\ 
3 & B_3 & 360 & 2 & 3\times A_4\\ \hline
4 & A_4 &  45 & 0 & A_4\\
4 & B_4 & 360 & 1 & A_4\\
4 & C_4 & 540 & 2 & 3\times A_4\\
4 & D_4 & 270 & 2 & 2\times A_4,A_6\\ \hline
5 & A_5 & 270 & 1 & A_4\\ 
5 & B_5 & 108 & 1 & A_6\\ 
5 & C_5 & 1080& 2 & 2\times A_4,A_6\\ \hline
6 & A_6 & 18  & 0 & A_6\\ 
6 & B_6 & 216 & 1 & A_6\\
6 & C_6 & 45  & 2 & 3\times A_4 \\ 
6 & D_6 & 270 & 2 & A_4, 2\times A_6 \\ 
6 & E_6 & 180 & 2 & 3\times A_6 \\ \hline
\end{array}$$
\end{table}
There is only one orbit for $k=2$. \qed

This gives also the information about the structure of the deep holes and
the cocode $K^6/{\cal C}_6$.

\begin{theorem}\label{deepholes}
The covering radius of the Hexacode ${\cal C}_6$ is $2$.
There is one type of deep holes in $K^6$. Representatives 
are the ${\rm Aut}({\cal C}_6)$-orbits $A_2$, $B_3$, $C_4$, $D_4$, $C_5$, $C_6$, $D_6$
and $E_6$. For every deep hole there are exactly three codewords with
distance~$2$. The three orbits $A_0$, $A_1$ and $C_6$
form a complete system of representatives for the  cocode $K^6/{\cal C}_6$, 
representing the cosets of minimal weight $0$, $1$ and $2$, respectively.
\end{theorem}
\an B:~\cite{CoSl-Orbit}; L: partially~\cite{CPS,BCQ,Bo-grass}; V: unknown.

The $135$ deep holes of weight $2$ are partioned into $45$ sets of ``trios'',
the members of each trio are representing the same coset in $K^6/{\cal C}_6$.
The subcode of ${\cal C}_6$ generated by pairs of members 
in a trio forms a frame which corresponds to a twisted construction of
${\cal C}_6$ from a $D_8^*/D_8$-code (cf.~the end of section~\ref{KBLV}).

\medskip

From the next theorem, one deduces immediately that
the $18$ vectors of weight $6$ in the Hexacode
are the smallest possible number of elements necessary to form a 
generalized $2$-design with $n=k=6$.
\begin{theorem}[Th.~5 and 6 in~\cite{BoCo}]\label{fischerun}
For the number of elements of a generalized \hbox{$2$-$(n,k,\lambda)$} design $Y$ 
of type~$3$ one has
$$ |Y|\geq \cases{3n, & for $k<n$,  \cr 2n+1,  & for $k=n$.}$$
\end{theorem}
\an B:~\cite{RaWi}, L:~\cite{DGS}, V: unknown.

The set of the $45$ weight $4$ vectors in the Hexacode has 
the smallest cardinality for a generalized $2-(6,4,\lambda)$ design.

By taking the $253$ of the $759$ vectors of weight $8$ in the binary Golay code 
having first coordinate~$1$, one gets 
the essentially only tight $4$-design~\cite{Br-elliptic}. 
The $196560$ vectors of squared length $4$ in the Leech lattice form
the only tight spherical $11$-design~\cite{BaDa1,BaDa2} in dimension greater 
then $2$. This leads to the question: 
Is there a good notion of tight generalized $t$-designs, using a
bound generalizing Theorem~\ref{fischerun} for its definition, 
characterizing one of the two designs belonging to the Hexacode? 

\section{Lexicographic codes}

The {\it lexicographic code\/} of length $n$ and minimal distance $d$ is 
defined by the greedy algorithm: After writing down the elements of $K^n$ in 
lexicographic order one chooses in every step the lexicographic first 
word which has distance at least $d$ to the already chosen codewords.

\begin{theorem}[Conway-Sloane~\cite{CS-lexi}]
The lexicographic code of length $2$ and minimal distance $2$ is $\epsilon_2$. 
The lexicographic code of length $6$ and minimal distance $4$ is
the Hexacode ${\cal C}_6$.
\end{theorem}
\an B:~\cite{CS-lexi}; L:~\cite{CS-lamin}.

Define {\it self-orthogonal lexicographic codes\/} by restricting the choice of 
the next codeword to the dual code of the code spanned by the codewords 
already chosen. This is some analogy to the definition of integral laminated lattices. 

\begin{theorem} 
The self-orthogonal lexicographic codes with minimal distance $1$, $2$,
$3$ and $4$ are ``periodic'' under direct sum. The periods are $1$, $2$,
$5$ and $6$ with periodicity elements 
$\gamma_1$, $\epsilon_2$, ${\cal C}_5$ and ${\cal C}_6$ respectively.
\end{theorem}
\an B: $c_2$, $H_8$, $g_{22}$ and $g_{24}$~\cite{Monroe}; 
L: $\bf Z$, $E_8$, $\Lambda_{23}$ and $\Lambda_{24}$~\cite{PP-intlamin,CS-intlamin}; 
V: $V_{\rm Fermi}$, $V_{E_8}$, $V\!B^{\natural}$ and $V^{\natural}$.

\section{Relations to binary codes, lattices and vertex operator algebras}\label{KBLV}

In this section, we assume that the reader is familiar with the notation 
of a vertex operator algebra (VOA) and a vertex operator super algebra
(SVOA) (see~\cite{FLM,FHL,Kac-VOA} for an introduction). All (S)VOA's are assumed
to be simple, unitary and ``nice'' (cf.~\cite{Ho-dr}, Ch.~1).

All the definitions and results of this work have analogies for binary codes, 
lattices and VOA's, although for VOA's the theory is not completely developed.
Analogously to the relation between binary codes and lattices and between 
lattices and VOA's one has two constructions (an ``untwisted'' and a 
``twisted'' one) for binary codes from Kleinian codes.

{\it Construction A:\/} Define a map $\rho_A$ from Kleinian codes of length $n$
to binary codes of length $4n$ by 
$$\rho_A(C):=\widehat{C}+d_4^n,$$
where $\widehat{\phantom{X}}:K^n \longrightarrow \F_2^{4n}$
is the map induced from $\hat{\phantom{.}}:
K\cong (D_4^*/D_4)\longrightarrow (D_2^*/D_2)^2\cong\F_2^4$,
$0\mapsto (0000)$, $a\mapsto (1100)$, $b\mapsto (1010)$, $c\mapsto (0110)$
and $d_4^n=\{(0000),(1111)\}^n$.
So every codeword in $C$ is replaced with $2^n$ binary codewords 
in $\F_2^{4n}$.

{\it Construction B:\/} Assume $n$ is even. Then
$$\rho_B(C):=\widehat{C}+(d_4^n)_0 \cup 
\widehat{C}+(d_4^n)_0+\cases{
(1000\ \ldots\ 1000\ 1000), & if $n\equiv0\pmod{4}$, \cr
(1000\ \ldots\ 1000\ 0111), & if $n\equiv2\pmod{4}$,}$$
where $\widehat{\phantom{X}}:K^n \longrightarrow \F_2^{4n}$ is the map 
as defined before 
and $(d_4^n)_0$ is the subcode of $d_4^n$ consisting of vectors of weight 
divisible by $8$.

\begin{lemma} If $C$ is a linear, self-dual resp.~even Kleinian code then
$\rho_A(C)$ is a linear, even self-dual resp.~doubly even binary code.
The same is true for $\rho_B(C)$ if the length is even. \qed
\end{lemma}
\begin{lemma}For the weight enumerators one has:
\begin{eqnarray*}
W_{\rho_A(C)}(x,y)& = & W_C(x^4+y^4,2x^2y^2), \\
W_{\rho_B(C)}(x,y)& = & \frac{1}{2} W_C(x^4+y^4,2x^2y^2)+\frac{1}{2}
(x^4-y^4)^n+\\ 
&& \qquad\frac{2^n}{2}\cdot \left((x^3y+xy^3)^n+(-1)^{n/2}(x^3y-xy^3)^n\right). 
\qquad \qed
\end{eqnarray*}
\end{lemma}
\an B--L: see~\cite{CoSl}, Ch.~7; L--V:~cf.~\cite{Ho-dr}, Ch.~1 and~5.

Remarks: 

$\rho_B(\hexa)$ gives the Golay code. (This is the MOG-construction.)

If we denote the untwisted (twisted) construction from binary codes to lattices
and from lattices to VOA's also with $\rho_A$ resp.~$\rho_B$ 
(cf.~\cite{DGH-virs}) then one has
$$\rho_{X}(\rho_{Y}(\rho_{Z}))=\rho_{\pi(X)}(\rho_{\pi(Y)}(\rho_{\pi(Z)})),
\quad\hbox{with $X$, $Y$, $Z \in \{A,B\}$ and $\pi\in S_3$.}$$

\smallskip

\noindent{\it Markings and frames:}

A {\it marking\/} for a code $C$ is the choice of a vector 
${\cal M}\in (K\setminus\{0\})^n$. 
Table~\ref{hexaorbits} shows that there exist $5$ inequivalent 
markings for the Hexacode.

For $i=1$, $\ldots$, $n$ we define
$$I_i=\cases{\{(4i-3,4i-2),(4i-1,4i)\}, & if ${\cal M}_i=a$, \cr
             \{(4i-3,4i-1),(4i-2,4i)\}, & if ${\cal M}_i=b$, \cr
             \{(4i-3,4i),(4i-1,4i-2)\}, & if ${\cal M}_i=c$.} $$
Then $I=\bigcup_{i=1}^n I_i$ is a {\it marking\/} for the binary code
$\rho_X(C)$ as defined in~\cite{DGH-virs}. As described in~\cite{DGH-virs}
one gets from $I$ a {\it $D_1$-frame\/} in $\rho_{X}(\rho_{Y}(C))$
(or equivalent a ${\bf Z}_4$-code, cf.~\cite{CS-z4}) and a 
{\it Virasoro frame} in $\rho_{X}(\rho_{Y}(\rho_{Z}(C)))$.
Since ${\rm Aut}(K^n)=S_3^n{:}S_n$ acts transitively on $(K\setminus\{0\})^n$
we can assume ${\cal M}=(aa\ldots a)$ by replacing $C$ with an equivalent code.
For this standard marking we define the {\it symmetrized (marked) weight
enumerator\/} ${\rm swe}_C$ as
$${\rm swe}_C(U,V,W)={\rm cwe}_C(U,V,W,W).$$

The {\it symmetrized marked weight enumerator\/} 
of the above marked binary code $\rho_X(C)$ as defined in~\cite{DGH-virs}
can be obtained from ${\rm swe}_C(U,V,W)$:
\begin{lemma}
\begin{eqnarray*}
{\rm smwe}_{\rho_A(C)}(x,y,z) & = & {\rm swe}_{C}(x^2+y^2,2xy,2z^2), \\
{\rm smwe}_{\rho_B(C)}(x,y,z) & = & \frac{1}{2}{\rm swe}_{C}(x^2+y^2,2xy,2z^2)+
\frac{1}{2}(x^2-y^2)^n+\\
&& \qquad\qquad \frac{1}{2}\cdot 2^n((x+y)^n+(-1)^{n/2}(x-y)^n)z^n. 
\qquad\quad \qed 
\end{eqnarray*}
\end{lemma}
\an B--L:~\cite{DGH-virs}; L--V:~\cite{DGH-virs}.

We remark that the symmetrized marked weight enumerator of an even self-dual
code belongs to a ring of polynomials with Molien series 
$\left(1+\lambda^4\right)/\left((1-\lambda^2)^2(1-\lambda^6)\right)$
generated by
$p_2={x^2} + 2\,{y^2} + {z^2}$, 
$q_2={x^2} + 4\,y\,z - {z^2}$,
$p_4={x^4} + 8\,{y^4} + 6\,{x^2}\,{z^2} + {z^4}$, 
$p_6= {x^6} + 6\,{x^2}\,{y^4} + 4\,{y^6} + 24\,{x^2}\,{y^3}\,z + 
   12\,{x^2}\,{y^2}\,{z^2} + 6\,{y^4}\,{z^2} + 8\,{y^3}\,{z^3} + 
   3\,{x^2}\,{z^4}$ 
subject to one relation for $p_4^2$.

\medskip 

Now, we describe how codes and lattices can be understood in terms of VOA's.
Let $V$ be a rational VOA whose intertwiner algebra is abelian, i.e.~the
set of irreducible $V$-modules form an abelian group $G$ under the fusion 
product (cf.~\cite{DoLe}).
The map $\alpha:G\longrightarrow \C^*$, $M\mapsto e^{2\pi i h(M)}$, where 
$h(M)$ is the conformal weight of the $V$-module $M$ defines a quadratic
form on $G$ and can be interpreted as an element of $H^4(K(G,2),\C^*)$; 
where $K(G,2)$ is the Eilenberg-MacLane space with $\pi_2(K(G,2))\cong G$
(see~\cite{Ho-simple}). Another description is the following:
The monodromy structure of the intertwiner operators 
of $V$ give rise to a three dimensional topological quantum field theory
which is example~I.1.7.2 of~\cite{Turaev}.

The fusion algebra of $V^{\otimes n}$ is ${\cal F}(V^{\otimes n})\cong {\bf Z}[G^n]$.
A subgroup $C\subset G^n$ is called an {\it even self-orthogonal linear code\/}
if $C$ is an isotropic subspace of the quadratic space $(G^n,\alpha^n)$.
It is proven in~\cite{Ho-simple} that (simple) VOA-extensions $W$ of 
$V^{\otimes n}$
are in one to one correspondence with such codes $C$; in particular, 
$W=\bigoplus_{\alpha\in C} M_{\alpha}$ has a unique VOA-structure up to 
isomorphism extending the VOA-structure of $V=M_0$.
The uniqueness follows from $H^3(K(C,2),\C^*)=0$.
Similar remarks hold for odd self-orthogonal codes and SVOA's. 

As an example, let
$V$ be the lattice-VOA $V_L$ belonging to an even integral positive definite
lattice $L$ of rank $n$. In this case $G=L^*/L$ with $\alpha$ induced from
$e^{2\pi i \frac{(.,.)}{2}}:\R^n\longrightarrow \C^*$,
where $(\,.\,,\,.\,)$ is the standard scalar product of $\R^n$. In fact,
the triple $(G,\alpha,n)$ is a complete invariant of the {\it genus\/} of $L$
(see~\cite{Ni-genus}).

Since the VOA belonging to the root lattice $D_4$ of ${\rm Spin}(8)$ has 
four irreducible modules with the conformal weights $0$ and three times $\frac{1}{2}$
and one has ${\cal F}(V_{D_4})\cong \Z[\Z_4]$, we get from the above example
the following description of Kleinian codes:

Even (odd) self-dual $K$-codes of length $n$ are the same as self-dual VOA's
(SVOA's) of rank $4n$ with sub-VOA $V_{D_4}^{\otimes n}$, the $n$-th tensor
product of the VOA associated to the Level-$1$-representation of the
affine Lie algebra ${\rm \widehat{Spin}}(8)$. The automorphism group of $K^n$
corresponds to the outer automorphism group of $V_{D_4^{\otimes n}}$ 
in the VOA-sense (Triality of ${\rm Spin}(8)$!); 
the group algebra $\Z[K^n]$ is the fusion algebra of $V_{D_4^{\otimes n}}$.

One has a similar description for binary codes in terms of the lattice-VOA
$V_{A_1}^{\otimes n}$. 

For $V$ be the (non rational) Heisenberg-VOA $V_h$ of rank $1$ on has
$G^n=\R^n$, $\alpha=e^{2\pi i \frac{(.,.)}{2}}$. Isotropic
subspaces are even integral lattices, i.e., we have a $1:1$-correspondence
between rank $n$ VOA's containing the Heisenberg-VOA $V_h^{\otimes n}
\cong V_{h^n}$ and even integral lattices.

The description of (marked/framed) Kleinian codes, binary codes and lattices  
in terms of VOA's is summarized in the next table. 

$$\begin{array}{l|c|c|c|c|c}
\mbox{Object} & \hbox{Rank} & \hbox{Sub-VOA (framed)} & \hbox{Group} & 
\mbox{Sub-VOA} & \hbox{Group}  \\ \hline
\mbox{$K$-codes} & 4n & V_{D_4,*}^{\otimes n} & 2^n{:}S_n & 
                  V_{D_4}^{\otimes n} &  S_3^n{:}S_n \\
\mbox{binary codes} & 2n & V_{D_2}^{\otimes n} & 2^n{:}S_n & 
                  V_{A_1}^{\otimes 2n} &  S_{2n} \\
\mbox{lattices} & n & V_{D_1}^{\otimes n} & 2^n{:}S_n & 
                  V_{h}^{\otimes n} &  {\rm SO}(n) \\
\mbox{VOA's} & n/2 & L_{1/2}(0)^{\otimes n} & \mbox{cf.~\cite{GH-stabil}} & 
                  {\rm Vir}_n &  ``{\rm Aut}({\cal F}({\rm Vir}_n))" 
\end{array}$$

\smallskip

Construction A (including marking/frames) can now be completely understood 
in terms of VOA's as indicated in following table of inclusions:

$$\begin{array}{ccc|ccc|ccc}
& \hbox{K--B} && & \hbox{B--L} & && \hbox{L--V} & \\ \hline
 V_{D_4,*}^{\otimes n} & \supset & V_{D_2}^{\otimes 2n}  &
 V_{D_2}^{\otimes n}   & \supset & V_{D_1}^{\otimes 2n} &
 V_{D_1}^{\otimes n}   & \supset & L_{1/2}(0)^{\otimes 2n} \\
 \cup & & \cup &  \cup & & \cup & \cup & & \cup \\
  V_{D_4}^{\otimes n}   & \supset & V_{A_1}^{\otimes 2n} &
  V_{A_1}^{\otimes 2n}  & \supset & V_{h}^{\otimes 2n} &
  V_{h}^{\otimes n}     & \supset & {\rm Vir}_n
\end{array}$$

\smallskip

For all four theories one has analogous basic objects. We display
their relations in Table~\ref{extremalobjects}.

\begin{table}\caption{Extremal odd Codes, Lattices and SVOA's}\label{extremalobjects}

$$ \begin{array}{|l|*{15}{c}|}\hline
\mbox{Rank}\phantom{\frac{|}{|}}
& \frac{1}{2} & 1 & \frac{3}{2} & 2 &\frac{5}{2} & 3 &\frac{7}{2} 
& 4 &  \frac{9}{2} & 5 & \frac{11}{2} & 6 &\frac{13}{2} & 7 &\frac{15}{2}  \\ \hline
\mbox{$K$-Codes} &   &   &   &   &   &   &   & \gamma_1 &   &   &   &  &  &  &  \\
             &   &   &   & &   &   &   & \da   &   &   &   &   &  &  &  \\
\mbox{$\F_2$-Codes} &   &   &   &c_2&   &   &   & c_2^2 &   &   &   &c_2^3&  &  &  \\
             &   &   &   &\da&   &   &   & \da   &   &   &   &\da  &  &  &  \\
\mbox{Lattices}&   &\Z &   &\Z^2&  &\Z^3&  &\Z^4   &   &\Z^5&  &\Z^6 & &\Z^7& \\
             &   &\da&   &\da &  &\da &  &\da    &   &\da&   &\da  &  &\da& \\
\mbox{SVOA's}& V_F & V_F^2 & V_F^3 & V_F^4  & V_F^5  & V_F^6 &
             V_F^7& V_F^8 & V_F^9 &V_F^{10}&V_F^{11}&V_F^{12}&
              V_F^{13}&V_F^{14}&V_F^{15}\\ \hline
\end{array}$$

$$ \begin{array}{|l|*{10}{c}|}\hline
\mbox{Rank}\phantom{\frac{|}{|}}
 & 8& 12 & 14 & 15 & \frac{31}{2} & 20 & 22 & 23 & \frac{47}{2} & 24 \\ \hline 
\mbox{$K$-Codes} & \epsilon_2 & \delta_{3}^+ & & &  & {\cal C}_5  &  & &  & {\cal C}_{6}\\
             &\da  &\da       &         & &   &  &      & &  &   \\
\mbox{$\F_2$-Codes} & e_8 & d_{12}^+ &(e_7+e_7)^+ & &   &  & g_{22} &  & & g_{24}\\
             &\da  &\da       & \da        & &   &  &      & &  &   \\
\mbox{Lattices} &E_8  &
          D_{12}^+&(E_7+E_7)^+& A_{15}^+& & & & O_{23} &  & \Lambda_{24} \\
            &\da & \da & \da & \da &   &   &  &   & &   \\
\mbox{SVOA's} &V_{E_8}&
        V_{D_{12}^+}&V_{(E_7+E_7)^+}&V_{A_{15}^+}&V_{E_{8,2}^+}& & &
&V\!B^{\natural}&V^{\natural} \\
\hline
\end{array}$$

The arrow $\da$ denotes construction $\rho_A$ and the rank of a Kleinian code of
length $n$ is defined as $4n$.
\end{table}

\medskip

{\it Final Remarks:\/} The way from Kleinian codes over binary codes and lattices 
to VOA's is not canonically given. There is no way to see what is the
next step. But in the other direction there is in some sense always a canonical 
choice: Consider the self-dual objects of rank $24$. There are always
two objects without ``roots'': An even and an odd one.\footnote{
In the case of vertex operator algebras the uniqueness
of the moonshine module $V^{\natural}$ and the odd moonshine module 
$V\!O^{\natural}$ is only a conjecture.}
Look at the even subobject 
of the odd one. Exactly one of its $4$ modules contains ``roots''. Take the 
direct sum of the even subobject and the ``root''-module and consider inside
the subobject generated by the ``roots''. It is a direct product of 
indecomposable objects. The next step is now represented by ``Codes'' over
the modules of one such indecomposable object. 

There is one more such step before Kleinian codes, namely codes over the 
$3$-dimensional topological quantum field theory belonging the
vertex operator algebra $V_{D_8}$. 

\bigskip

\small

\noindent{\it Some historical comments and further developments:\/}

I found the structure of Kleinian codes as developed in this paper
by searching for an analogue of the shorter Moonshine module in autumn 1995.
This was motivated by the work on Virasoro frames inside the Moonshine module.
The weight enumerator of the shorter Hexacode (which is not a ${\bf F}_4$-code)
dropped out. Compare the last paragraphs above. 

A first outline of this paper was distributed during the first two month  
of 1996 including all the results but most proofs not yet written up in Kleinian 
code language. Some other preliminary versions, but now without Section~5, were 
distributed in summer 1996. The only exception to this is the extremal code of 
length~12. I tried to find such a code by hand (cf.~letter to
Hirzebruch~\cite{Ho-hiletter}),
but without success. Back in Germany in October 1996, it popped up on the 
screen of my old \hbox{AT-286 PC} 
after a few minutes (or hours) by running a simple back-tracking algorithm.
This code was also found in~\cite{CRSS-quant}, where the authors applied 
the theory of Kleinian codes to quantum codes. This paper became the 
stimulus of a lot of research on quantum codes.
It seems that only a late 1996 preprint found the widest
distribution. I am sorry about the delay in publishing the paper.
I like to thank C.~Bachoc, J.-L.~Kim and V.~Pless for comments
on the final version.

\medskip

Since that time,  Kleinian codes have been investigated further.
In the following, I will give an overview. 

\smallskip

{\it Section 2:\/} The invariant ring for the complete weight enumerator
of even self-dual Kleinian codes has been given in~\cite{RaSl-hand}. 

\smallskip

{\it Section 3:\/} Examples of cyclic self-dual codes for all odd length 
have been given by M.~Ran and J.~Snyders in~\cite{RaSn}.

It was pointed out to me by J.-L.~Kim that the 
papers~\cite{GKHP-extremal2,BaGa-extremal} are answering partially my question for 
the smallest codes with trivial automorphism group: There is at least one 
such code of length~$12$ (called $QC\_12g$ in~\cite{GKHP-extremal2}; non even) 
and there are at least $273$ such extremal even codes of length~$14$ 
(see~\cite{BaGa-extremal}). Since all the even codes of length~$8$ and $10$
without weight~$2$ vectors are extremal, it follows from 
Section~3 and~\cite{BaGa-extremal} that the answer 
for even self-dual codes must be $12$ or $14$.

\smallskip

{\it Section 4:\/} The upper bound of Theorem~\ref{bound-odd} has been
sharpened by E.~Rains in~\cite{Ra-shadow} to $d\leq 2[n/6]+2+e$ with 
$e=1$ for $n\equiv 5 \!\!\pmod{6}$ and $e=0$ else. For $6|n$, a code meeting this 
bound is even. An analogue sharpened bound for binary codes can also be 
found in~\cite{Ra-shadow} and for odd lattices in~\cite{RaSl-lattice}.

In~\cite{GKHP-extremal,GKHP-extremal2}, Gaborit, Huffman, Kim and Pless 
classified self-dual Kleinian codes with minimal weight reaching the above
bound for length $8$, $9$ and $11$ (there are $5$, $8$ resp.~$1$ such codes). 
They also proved the uniqueness of the extremal even code of length~$12$.
There is no such code for length $13$ (see~\cite{RaSl-hand}).
For even codes, the length $10$ has been settled in~\cite{BaGa-extremal} ($19$ codes), 
where also partial results for length~$14$ and $18$ are obtained.  

{\it Section 5:\/} C.~Bachoc (see~\cite{Ba-nonbin}) has proven 
Theorem~\ref{design} and some extensions of it for all $n$ by using discrete
harmonic analysis on $X_k$. 
Interestingly, this approach works only for alphabets with $2$, $3$ and 
$4$ elements and a unique choice of group structure and bilinear form. 
For four elements, one gets our scalar product on $K$. The binary analogue was
studied before in~\cite{Ba-harm}. This approach forms the direct analogue 
to the approach of B.~Venkov for lattices~\cite{Ve-design}.


\providecommand{\bysame}{\leavevmode\hbox to3em{\hrulefill}\thinspace}

\end{document}